\documentclass[10pt,twoside]{article}
\usepackage{amssymb}
\usepackage{latexsym}
\usepackage{psfig,latexsym,rotating,amsmath,epsfig}
\usepackage{here}
\def\N{\mathbb{N}}
\def\Z{\mathbb{Z}}

\def\C{\mathbb{C}}

\def\Sl{\mathfrak{S}}
\def\E{\mathfrak{E}}

\def\1{\mathbf{1}}
\def\:{\lrcorner}
\def\#{\sharp}

\def\i{\iota}
\def\l{\lambda}
\def\a{\alpha}
\def\b{\beta}

\def\r{\rho}
\def\o{\circ}

\def\x{\otimes}
\def\<#1,#2>{\langle#1,\,#2\rangle}

\def\K{\mathbb{K}\,}
\def\S{\mathbb{S}\,}

\def\qed{\ensuremath{\quad\Box\quad}}
\def\pfill{\par\vskip2mm plus1mm minus1mm\noindent}
\def\inv#1{\raise.1em\hbox to 0pt{$^{-1}$\hss}_{#1}\;}

\newtheorem{Theorem}{Theorem}
\newtheorem{Lemma}[Theorem]{Lemma}
\newtheorem{Corollary}[Theorem]{Corollary}
\newtheorem{Definition}[Theorem]{Definition}
\newtheorem{Beispiel}[Theorem]{Example}

\begin{document}

\setlength{\textwidth}{5in} \setlength{\textheight}{7.8in}

\setlength{\baselineskip}{13pt}

\date{}

\title{INVARIANTS OF 3-MANIFOLDS FROM REPRESENTATIONS OF THE FRAMED-TANGLE CATEGORY}

\author{OLAF M\"ULLER\\Max-Planck-Institut f\"ur Mathematik in den Naturwissenschaften\\Inselstrasse 22-26, 04103 Leipzig, Germany\\Olaf.Mueller@mis.mpg.de}
\maketitle
\begin{center}Received          \\  \end{center}

\bigskip

\begin{abstract}  
\noindent We will construct a monoidal functor ("a monoidal representation") from the category of framed tangles into the tensor category  over a fixed ground vector space which is invariant under Kirby moves and so gives rise to an invariant of 3-manifolds.
\end{abstract}

\section{Introduction}

\begin{figure}[H]
\epsfig{file=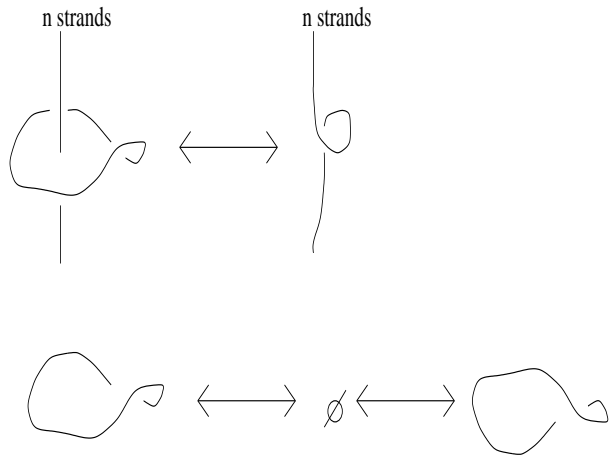,height=6cm,width=8cm,clip=,angle=0}
\caption{\footnotesize{Kirby moves (Fenn-Rourke moves). Two framed links represent one and the same 3-manifold iff you can link them by a finite sequence of these moves. The n strands are always thought to be parallel}}
\label{durchbruch}
\end{figure}

It is a well-known fact due to Kirby ~\cite{ki} (most comprehensible proof: ~\cite{nl} ) that 3-manifolds are in 1:1-correspondence to framed links quotiented by a set of combinatorial moves, the so-called Kirby moves. In this situation, it would be a quite obvious idea to subdivide a standard projection of the framed link by a quadratic grid so that in each small square the link is isotopic to some standard framed tangle and then to represent the tangles by morphisms in the tensor category over a fixed ground vector space $V$. Each strand is represented by a copy of $V$, the composition in the framed tangle category $\Sl'$ becomes the composition of linear maps, the juxtaposition of tangles becomes the tensor product of maps; the functor is therefore called a {\em monoidal} representation. We will give a sufficient and easily verifiable condition on a monoidal representation of $\Sl'$ to factorize through the Kirby moves and to generate therefore an invariant of 3-manifolds. A priori this is very difficult to check because in the local version of the Kirby moves due to Fenn and Rourke ~\cite{fr} (cf. fig. \ref{durchbruch} ) there are infinitely many moves parametrized by the number of non-closed strands giving rise to infinitely many conditions on the generating morphism of the representation, the so-called S-matrix representing a single crossing.

\noindent Other approaches to this problem are due to Wenzel ~\cite{we} and Turaev ~\cite{t3} and deal with colourings of framed tangles giving up the possibility to construct one single representation and to gain skein relations that already turned out to be very useful in knot theory.  

\bigskip

\noindent Wenzel defined the Markov traces, a family of Markov-invariant functionals on $ B_\infty $ . The twists are encoded by modelling each component as a ribbon of several strands; the number of strands is called the {\em colouring} of the component, and to calculate the invariant one has to sum up combinations of colourings. Therefore the invariant has no representation properties: Consider a single crossing where the upper strand is coloured red, the lower one green, then the square of the coloured braids does not emerge in the sum of colourings of the squared braid because the colourings do not fit together. 

\bigskip

\noindent In Turaev's approach  by modular categories, going out from monoidal representation of coloured tangles, there is the same problem, and a new difficulty arises: The objects (strands) are coloured by simple objects.  

\noindent Consider modules over ground ring $\K$ . An object $V$ is called {\em simple} iff the map  
$k \mapsto k \otimes id_V $

\noindent is a bijection from $ \K $ to $ End(V) $ . 

\noindent This implies that if $ \K $ is a field then $ V \cong \K $ as fields and each tangle $  L $ is represented by an endomorphism of the ground field; therefore $ L \otimes id = id \otimes L $ so that it is no more possible to localize strands. That means that over a field Turaev's construction does not work.

\bigskip

\noindent At this point the simplest idea would be to construct monoidal representation of the category of framed tangles that factorize through the Kirby moves. For this goal it is necessary to reduce the infinitely many Fenn-Rourke conditions to a finite set of conditions. This is exactly what will be done in section 8, all sections before just recall more or less well-known facts which can be found in the given references.  
 
\bigskip

\noindent One main result which will be presented in section 9 is that {\em irreducible} representations of $\Sl '$ that satisfy only the two first pairs of Fenn-Rourke conditions generate invariants of 3-manifolds. But the proof of the existence of such a representation still fails at the moment because there is still no classification of solutions of the Yang-Baxter equation for $ dim(V) \geq 3$.

\bigskip

\noindent A part of this paper is part of my diploma thesis at the University of Bonn, and I would like to thank Prof. C.F. B\"odigheimer who was my advisor in that time.

\section{From 3-Manifolds to Framed Tangles}

\noindent The concept of Dehn surgery became possible by the technique of Heegard decomposition on the one hand and by the existence of a simple set of generators of the surface mapping class groups on the other hand. Since we know that we can get any 3-manifold by gluing two handlebodies along a homeomorphism of their boundaries and that the topological type of the resulting 3-manifold only depends on the {\em isotopy} class of the surface homeomorphism, the mapping class group of a surface becomes interesting. A generating set are the {\em Dehn twists} along curves on the surface. So we can decompose the homeomorpism and study the effect of a single Dehn twist that consists in eliminating and re-gluing a torus from $\S^3$ 

\begin{Theorem}
Any compact, orientable 3-manifold can be produced from $\S^3$ by eliminating of some embedded tori and re-inserting them by gluing along some homeomorphisms of the boundaries (Dehn surgery).
\end{Theorem}

\noindent {\em Proof} cf. ~\cite{ps} .

\bigskip

\noindent The isotopy class $[f]$ of a homeomorphism $f$ is encoded by the $2 \times 2$-matrix of $f_*$ in the fundamental group $ \Z \oplus \Z $: 

\begin{displaymath}
\mathbf{R} = 
\left( \begin{array}{cc}
1 & 0  \\
p & 1  
\end{array} \right)
\end{displaymath}
    
\noindent where $ p \in \Z $. Therefore $[f]$ is determined by $p \in \Z $. 

\noindent Therefore you can encode $[f]$ by the core $l$ of the torus and the winding number $p$ of the meridian, i.e. by embeddings of ribbons (by a {\em framed} link)

\begin{Theorem}

Two framed links represent one and the same 3-manifold iff you can link them by a finite sequence of the moves in fig. \ref{durchbruch}.

\end{Theorem}

\noindent {\em Proof.} Difficult result of Kirby ~\cite{ki} , improved by Fenn and Rourke ~\cite{fr} , the most comprehensible proof is due to Ning Lu ~\cite{nl} \qed 

\bigskip

\noindent We cite some known results that can be found in ~\cite{ps}, for example. 

\begin{Theorem}
Let $M$ resp. $ M'$ be the 3-manifolds associated to the framed links $L$ resp. $L'$, then $ L \sqcup L' $ induces the manifold $ M \# M' $ \qed  
\end{Theorem}

\begin{Beispiel}
The trivial link of m components with trivial framing (no twist) induces the manifold $ ( \S^1 \times \S^2 ) ^{ \# m} $.
\end{Beispiel} 

\begin{figure}[H]
\epsfig{file=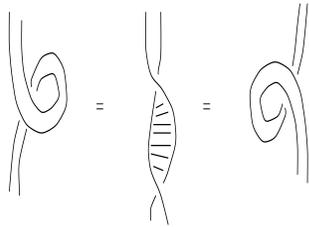,height=3cm,width=4cm,clip=,angle=0}
\caption{\footnotesize{Framings of a vertical strand. The back side is hatched}
\label{rahmungen}}
\end{figure}

\noindent Given a framed link, we now choose a {\em regular projection} of it, i.e. an embedding of the corresponding union of ribbons representing this framed link whose heighth extrema in the projection plane are isolated points and that is "nearly parallel" to a projection plane, what means that the normal vector of the ribbons differs only slightly from the one of the projection plane. So the ribbons turn always onre and the same side to the observer. Figure \ref{rahmungen} illustrates that any framed link has such a projection. By this convention we can represent the framed link only by its core.

\begin{figure}[H]
\epsfig{file=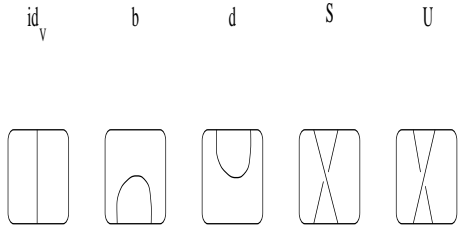,height=3cm,width=6cm,clip=,angle=0}
\caption{\footnotesize{Monoidal generators for $ \Sl $ and corresponding linear maps}}
\label{rasterfahndung}
\end{figure}

\noindent Now, we subdivide such a projection by a quadratic lattice in such a way that the part of the framed link within a small sqare of the lattice is isotopic to one of the five standard framed tangles in fig. \ref{rasterfahndung} that generate by composition and juxtaposition the whole category of framed tangles. Then we replace them following the scheme in the figure where each strand is represented by one copy of $V$, translating composition and juxtaposition to composition and tensor product of linear maps, $S$ and $U$ being linear endomorphisms of $V \x V$, while $b$ and $d$ representing caps and cups are linear maps $ \K \rightarrow V \x V $ and $ V \x V \rightarrow \K$, respectively. For a general quadruple $(S,U,b,d)$ this assignment, called $\rho_{S,U,b,d}$ from now on, is only well-defined on the level of projections. The next theorem will give conditions on it to be still well-defined for framed links.

\noindent So in the end, we have got to care about the category of framed tangles. For technical details about this cf. ~\cite{t2} and ~\cite{ps}, for example. As the framed tangles form a {\em monoidal category} we will consider this categorial concept now.

\section{Monoidal Categories and Monoidal Representations of Framed Tangles}

A {\em (strictly) monoidal category} is a category $ K $ endowed with a covariant associative functor $ \x : K \times K \rightarrow K $ with unity. Covariant means

$ (f \o f' ) \x (g \o g' ) = ( f \x g ) \o ( f' \x g' ) $,

$ id_V \x id_W = id_{V \x W} $.

\noindent Associativity means on the level of objects that for each three objects $U,V,W$ :

$ (U \x V) \x W = U \x ( V \x W)$

and on the level of morphisms that for each three morphisms $f,g,h$ :

$ (f \x g) \x h = f \x ( g \x h) $.

\noindent There is a unity, i.e. an object $ \1 $ with:

$ V \x \1 = V = \1 \x V $ for all objects $V$ and

$ f \x id_{ \1 } = f = id_{ \1 } \x f $ for all morphisms $f$.

\bigskip

Examples are 

\begin{itemize}

\item{the tensor category $\E _V$ over a fixed ground vector space $V$,}
\item{the category of tangles $ \Sl $,} 
\item{the category of framed tangles $ \Sl ' $.}
\

\end{itemize}

\noindent A functor between two monoidal categories is called {\em monoidal} iff it translates the tensor product of the first category into the one of the second category. An example is the obvious forgetful functor $ \Sl'  \rightarrow \Sl $. Monoidal functors with image in $ \E_V $ for some $V$ we will call {\em monoidal representations}.
From now on, let $V$ be a finite dimensional vector space over a field $ \K $, $v:= dim_{\K}(V)$. The standard example will be $ \K = \C $.

\bigskip 

\noindent For the following, it will be convenient to have another definition at hand, namely the one of a {\em partial trace} of an endomorphism in $ \E_V $. 

\noindent If we consider the effect of the closure of the last strand in the representation it turns out that it is just taking the sum over the last tensor index, similar for the first strand. So, chosen a basis $B = e_1...e_m$ of $V$, we define

$ r_B ,l_B : End(V^{\otimes i}) \rightarrow End (V^{\otimes i-1})$

 and for $A \in End(V^{\otimes i})$

$$r_B(A) (e_{j_1}...e_{j_{i-1}}) := \sum_{k_1...k_{i-1},r_i = 1}^m A_{j_1...j_{i-1}r_i}^{k_1...k_{i-1}r_i}    
e_{k_1} \otimes ... \otimes e_{k_{i-1}} \qquad ,$$ 

$$ l_B(A) (e_{j_1}...e_{j_{i-1}}) := \sum_{l_1,k_1...k_{i-1}=1}^m A_{l_1j_1...j_{i-1}}^{l_1k_1...k_{i-1}}    
e_{k_1} \otimes ... \otimes e_{k_{i-1}} \qquad .$$ 

\noindent $r_B$ and $l_B$ are called right and left partial trace, respectively, and are in some way a n-th root of the trace.

\noindent As it stands, this definition depends on the choice of the basis $B$. For technical details as $Gl(V) \times Gl(V)$-covariance of partial traces cf. ~\cite{m} , for example.

\begin{Theorem}
\label{darstellung}
 Given linear maps

$ S : V \x V \rightarrow V \x V $ invertible

$ b : \C \rightarrow V \x V $

$ d : V \x V \rightarrow \C $. 

If and only if

\bigskip

\noindent (i) $S$ is a solution of the Yang-Baxter equation (YBE): 

$  (S \otimes id_V ) \circ (id_V \otimes S) \circ (S \otimes id_V) = (id_V \otimes S) \circ (S \otimes id_V ) \circ (id_V \otimes S) $

\noindent (ii) $d$ is symmetrical and has an orthonormal basis $ B= \{ e_i \} $; for this basis holds:

$ b ( \l ) = \l \sum_i e_i \x e_i $

\noindent (iii) $ (S^{ \pm 1 } \x id_v ) \o ( id_v \x b ) = ( id_v \x S^{ \mp 1 } ) \o ( b \x id_v ) $ 

\noindent (iv) $r_B(S) = l_B(S)$

\bigskip
\noindent then $ \rho_{S,S^{-1},b,d}$ is invariant under isotopies of framed tangles and therefore a monoidal functor $ \Sl ' \rightarrow \E_V  $ (a monoidal representation of $ \Sl '$). In this case we call $S$ an S-matrix.
\end{Theorem}

The conditions above are illustrated in the following figure.
      
\begin{figure}[H]
\epsfig{file=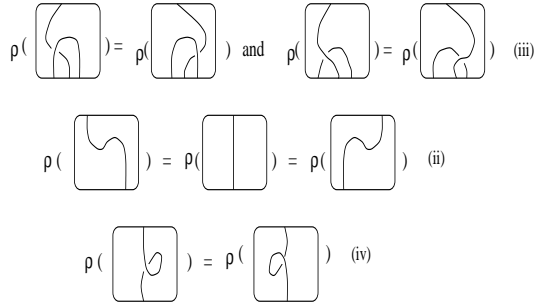,height=4cm,width=7cm,clip=,angle=0}
\caption{\footnotesize{Conditions on $ S,T,b $ and $ b,d $}}
\label{gleiten}
\end{figure}

\begin{figure}[H]
\epsfig{file=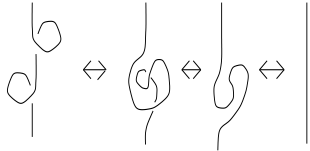,height=2cm,width=4cm,clip=,angle=0}
\caption{\footnotesize{Whitney trick}}
\label{wt}
\end{figure}   

\noindent {\em Proof.} It must be checked that any isotopy of framed tangles can be modelled by the moves corresponding these relations. But an isotopy of a ribbon is in particular an isotopy of its core. So you can use the corresponding result about representations of the (non-framed) tangle category found by Turaev ~\cite{t2}. Each time when there is a move corresponding to a Reidemeister-I-move replace it by the move corresponding to the Whitney trick using only Reidemeister-II and Reidemeister-III moves (cf. fig. \ref{wt}), let one loop be so big as you want it to be and the other one very small contained within a normal neighborhood of the ribbon, then stop the isotopy for a moment and lead the loop always within this neighborhood along the ribbon to a fixed point of the isotopy where you collect all small loops in this way. Finally, the Whitney trick reduces the number of loop types to two, and the only remaining move is the one in fig. \ref{gleiten}, lowest row. And this move is equivalent to the last condition:

\noindent Denote a (n,n)-tangle by $P$, the tangle we get if we close the last strand of $P$ by $P'$, and the one we get if we close the first strand by $'P$. Then 

$\rho_{S,S^{-1},b,d} (P') = r_B( \rho_{S,S^{-1},b,d}(P))$ , 

$\rho_{S,S^{-1},b,d} ('P) = l_B( \rho_{S,S^{-1},b,d}(P))$ 

\noindent because in terms of $b$ and $d$ with the definitions

$b^{(n)} := id_{n-1} \otimes b$ , $d^{(n)} := id_{n-1} \otimes d$ , 

$b^{(1)} := b \x id_{n-1} $ , $d^{(1)} := d \x id_{n-1} $ , 

\noindent one can write $r(A), l(A)$ as

\begin{figure}[ht]
\epsfig{file=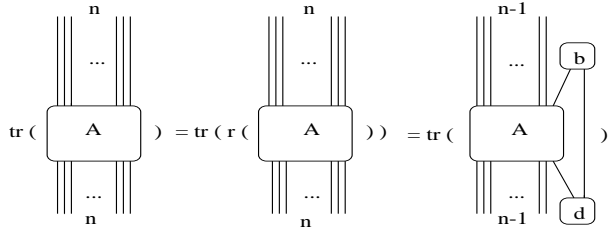,height=3cm,width=8cm,clip=,angle=0}
\caption{\footnotesize{An n-th root of the trace}}
\label{r}
\end{figure}

$r_B(A) = b^{(n)} \circ (A \otimes id) \circ d^{(n)}$ ,

$l_B(A) = b^{(1)} \circ (id \x A) \circ d^{(1)}$ (cf.  fig. \ref{r} )  \qed

\noindent As a consequence of Theorem \ref{darstellung} (ii) the choice of $b$ and $d$ is uniqely determined by the choice of a basis, and if $S$ is a matrix, then the notation $ \rho_S $ is well-defined.

\section{A symmetry of S-Matrices}

The natural scalar product $\langle .,. \rangle _n $ on $ V^{ \x n } $ is defined by

$ \langle v_1 \x ... \x v_n , w_1 \x ... \x w_n \rangle _n := d( v_1 , w_1 ) \cdot  ...  \cdot d( v_n , w_n ) $ .

\noindent Then we have

\begin{figure}[h]
\epsfig{file=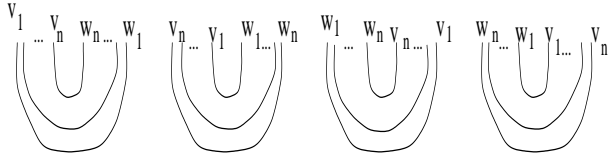,height=2cm,width=8cm,clip=,angle=0}
\caption{\footnotesize{Four graphs that represent the n-th scalar product}}
\label{skalar}
\end{figure}

$ \langle .,. \rangle _n = \rho (U_n) \o  (id_n \x T_n) = \rho (U_n) \o ( T_n \x id_n) $,

\noindent where $U_n$ is the tangle in fig. \ref{skalar}, $T_n$ the twist in the tensor category converting the order in each tensor product.

\begin{Theorem}
Let $ s $ be the rotation with angle $ \pi $ in the projection plane, 
$F$ an arbitrary (n,m)-tangle. Then for any monoidal representation $ \rho $ holds

$ ( T_m \o \rho (F) \o T_n )^+ = \rho ( s (F) ) $

\end{Theorem}

\noindent {\em Proof.\/} $ \langle \rho (F) (v) , w \rangle _m 
= \rho (B) ( v \x T_m (w) ) 
= \langle T_n v , \rho ( s (F) ) (T_m w) \rangle $

$= \langle T \o \rho ( s (F) ) (T_m w) , v \rangle $ \qed

\begin{figure}[h]
\epsfig{file=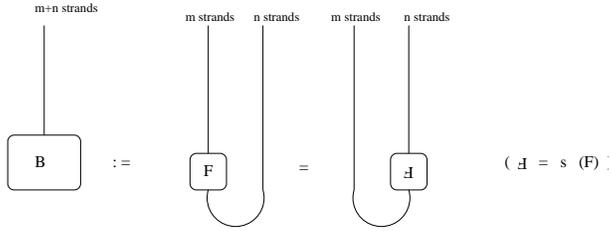,height=3cm,width=8cm,clip=,angle=0}
\caption{\footnotesize{Scalar product and rotation of tangles}}
\label{rumdrehen}
\end{figure}

\noindent As $ T^2 = \1 $ and $T$ is symmetrical, with the definition

\noindent $ A \in End(V^{ \x n })  $ T-symmetrical $ : \Leftrightarrow  T_n \o A $ symmetrical we get

\begin{Corollary}

Let $H$ be a $s$-symmetrical tangle, then $ \rho (H) $ is T-symmetrical; in particular $S$ is T-symmetrical  \qed

\end{Corollary}

\noindent In Coefficients this reads: $ (S^{ \pm 1} )_{ab}^{cd} = (S^{ \pm 1})_{dc}^{ba} $
(" Index rotation " )

\noindent Now we have already (cf. fig. \ref{gleiten} and Theorem \ref{darstellung} (iii))

$ S_{kl}^{pm} = \tilde S_{lm}^{kp} $ und $ \tilde S_{kl}^{pm} = S_{lm}^{kp} $

\noindent So instead of the two demands in \ref{darstellung} (iii) we can impose equivalently only one sign chose and the index rotation (which is a {\em linear} condition). Later, T-symmetry will be used in the proof of Theorem \ref{sym}.

\section{S-Matrices and Symmetries of the YBE}

It is a known fact that from a solution of the YBE you can derive other solutions:

\begin{Theorem}
\label{verwursten}
If $ S : V^ {\otimes 2} \rightarrow V^ {\otimes 2} $ is a solution of the YBE, then also

(i) $ \l S, \l \in \C $

(ii) $ S^{-1} $

(iii) $S^+$ (the transposed of $ S $ )

(iv) $ T \circ S \circ T $

(v) $ (Q \otimes Q) \circ S \circ (Q^{-1} \otimes Q^{-1}) $, 

where $ Q $ is an arbitrary invertible endomorphism $ V \rightarrow V $ .
\end{Theorem} 

\bigskip

\noindent Now the T-symmetry found in the previous section sheds some light on how S-matrices react to the transformations in Theorem \ref{verwursten}:

\begin{itemize}

\item{(because of T-symmetry):  Transposition = Conjugation with $T$} 
\item{(if $S$ yields a 3-manifold invariant $\iota_S$):  $ \rho_{S^{-1}} (M) = \rho_S (-M) $, 

where $-M$ is the manifold $M$ with orientation inverted}
\item{Scalar multiplication: Assume that $S$ and $ \lambda S $ are both S-matrices, then

$ \l S_{ab}^{cd} = \l^{-1} \tilde S_{bd}^{ac} = \l^{-1} S_{ab}^{cd} $

(second equation: index rotation for $S$ , first equation: index rotation for $ \l S $), 

and therefore $ \l = \pm 1 $, i.e. S-matrices are uniquely scaled up to sign.}
 
\end{itemize}

\section{Examples of Monoidal Representations for $\qquad$ $dim_{\C}(V) = 2$ }
 
We use the classification of solutions of the YBE in $dim(V)=2$ by Doll ~\cite{do} .

\bigskip

\noindent It turns out (cf. ~\cite{m}) that the only representation of that kind not satisfying the trivializing relation $ S^2 = \1 $ is generated (in the sense of \ref{verwursten} ) by

\begin{displaymath}
\mathbf{S} =  
\left( \begin{array}{cccc}
 k & 0 & 0 & 0 \\
 0 & 0 & q & 0 \\
 0 & p & 0 & 0 \\
 0 & 0 & 0 & k
\end{array} \right)
\end{displaymath}

where    $ k^2 pq = 1 $.

\bigskip

\noindent But in this family of solutions there is no invariance under Kirby moves.

\section{Gaining Skein Relations}

\begin{figure}[H]
\epsfig{file=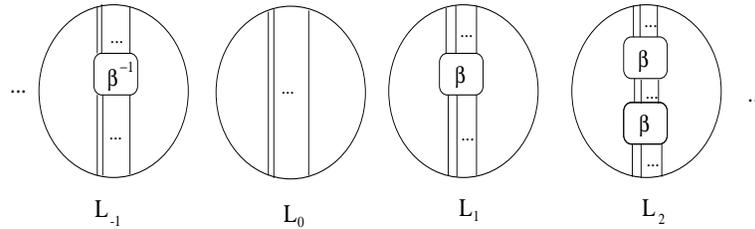,height=3cm,width=10cm,clip=,angle=0}
\caption{\footnotesize{A $ \b $-sequence is a sequence of framed links which have regular projections that are identical outside a circle and differ within by potences of a framed braid.} }
\label{bseq}
\end{figure}

\bigskip
\noindent Now we describe the usual procedure to get richer invariants from the above ones.
The YBE is homogenous of degree 3 in the entries of $S$, all additional conditions are of degree $ \leq 2 $, therefore the family of S-matrices forms a subvariety of $Mat(n^2 \times n^2 , \K)$ that can be parametrized locally by a polynomial as $S(x_1,...,x_n)$. So you can define an invariant $ \rho $ with values in the polynomial ring over the ground ring by setting 

\bigskip

$ \rho (M) (x_1,...,x_n) := \rho_{S(x_1,...,x_n)} (M) $.

\bigskip

\noindent The following theorem describes how you can gain skein relations from algebraic relations of the S-matrix. In the case that all S-matrices of a subvariety as above satisfy these algebraic solutions the skein relations can be transferred to the polynomial invariant associated to this subvariety.

\begin{Theorem}
\label{algeo} 
Given a solution $ S $ a S-matrix, $ \b $ a framed braid. For $ A:= \r _S ( \b ) $, any polynomial relation $ P(A) = 0 $, where  
$ \displaystyle{ P = \sum_{i=m}^n a_i X^i  \in \K [X] } $ (Laurent ring), implies

$$ \sum_{i=m}^n a_i \rho_S (L_i)  = 0 \qquad , $$  

where $L_i$ is an arbitrary $ \b $-sequence  \qed

\end{Theorem}   

\bigskip

\noindent A simple example are Conway sequences where the number of strands is 2 and therefore only potences of the unique generator $S$ itself appear. For a fixed  $ \b $ there is exactly one generating relation, namely the one given by the minimal polynomial of $ A = \r _S ( \b ) $ , that is in the case of Conway sequences the minimal polynomial of $S$.

\section{Invariance under Fenn-Rourke Moves}

From now on, we denote the upper right tangle in fig. \ref{durchbruch} by  $tw_n$; $tw_1=:tw$.

\noindent Define  $ S^{(i,j)} := id_i \x S \x id_j $, $tw^{(i,j)} := id_i \x tw \x id_j$, for $b$ and $d$ analogously. 

\begin{figure}[H]
\epsfig{file=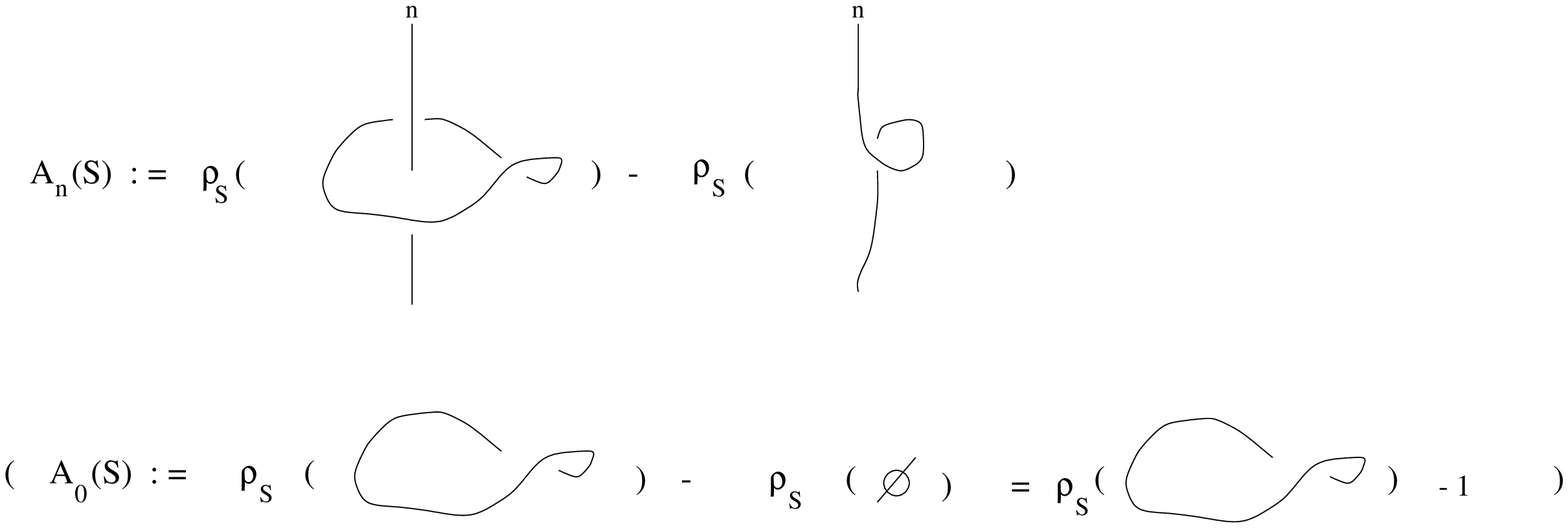,height=4cm,width=10cm,clip=,angle=0}
\caption{\footnotesize{Conditions from Fenn-Rourke moves}}
\label{AnBn}
\end{figure}

\noindent Define maps $A_n, B_n : V^{ \x n } \rightarrow V^{ \x n }  , n \in \N_0 $, 

\noindent $A_n$ as in fig. \ref{AnBn}, the maps $B_n$ analogously as mirror images.

\noindent Because of the Theorem of Fenn-Rourke (cf. ~\cite{ps}) we have to look for a matrix $S$ with $ A_n(S) = 0 \qquad \forall n \in \N_0 , B_1(S) = 0 $. 

\bigskip

\begin{figure}[H]
\epsfig{file=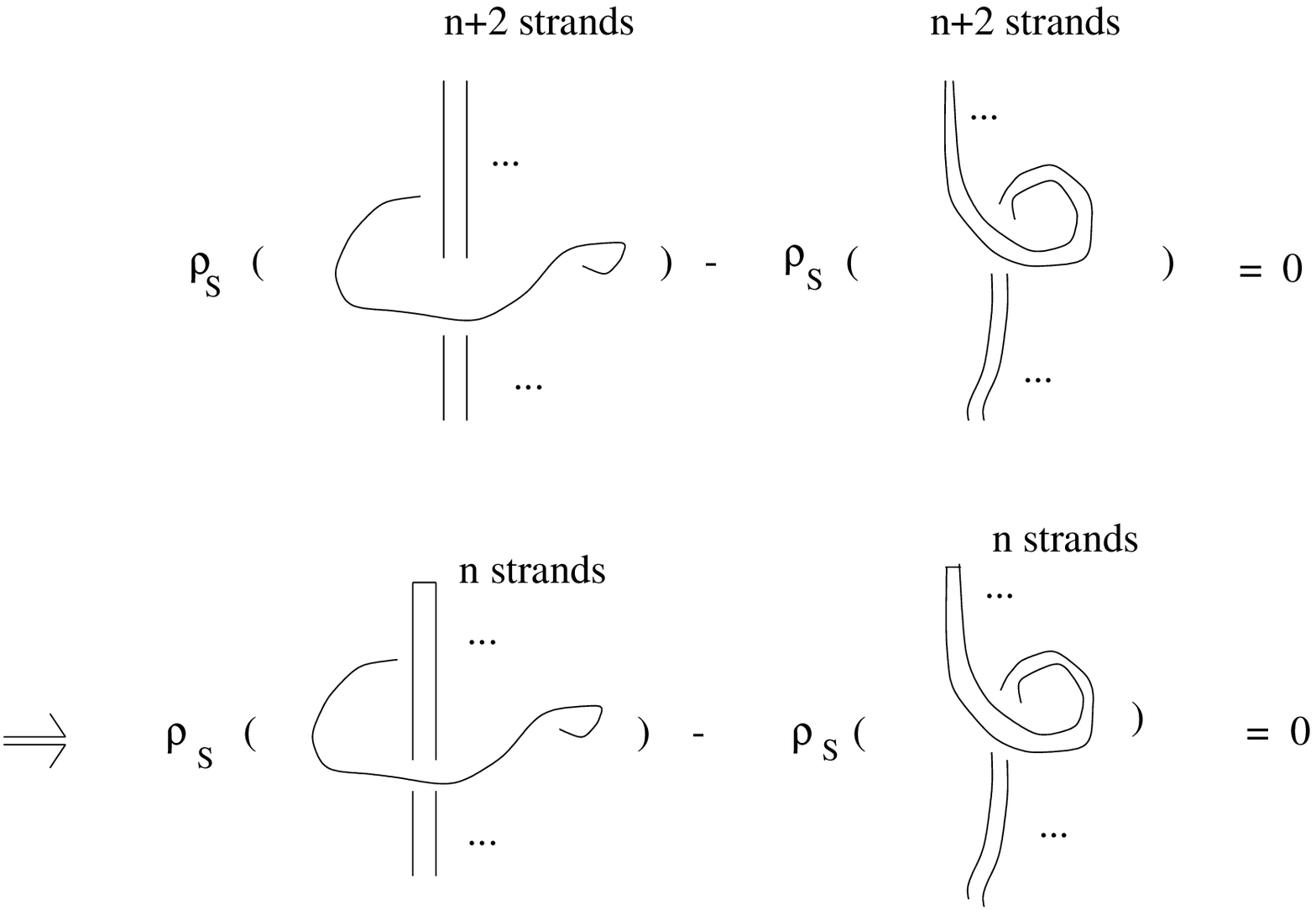,height=5cm,width=7cm,clip=,angle=0}
\caption{\footnotesize{$A_{n+2} = 0$ implies $A_n = 0$ }}
\label{durchziehen}
\end{figure}

\noindent It is easy to see that $ A_{n+2} (S) = 0 $ implies $ A_n (S) = 0$. All you have to do is to consider $ A_{n+2} (S) \o b^{(i)}  $ (cf. fig. \ref{durchziehen} )  and to pull the corresponding pair of closed strands down through the tangle. This isotopy does not change the resulting value, namely zero, of the representation. On the other hand, $b$ has full rang, so you are done.   

\noindent Now we give the main theorem concerning invariance under Fenn-Rourke moves.

\bigskip
\begin{Definition}

An endomorphism $E : V^{ \x n+2} \rightarrow V^{ \x n+2}$ with
      
$\qquad$       ($\a$) $ b^{(i,n-i)} \o E = 0 = E \o d^{(i,n-i)} \qquad \forall 0<i<n $
   
$\qquad$       ($\b$) $r(E) = 0 $
       
$\qquad$   ($\gamma$) $E$ commutes with $ S^{(i,n-i)} $ and $ tw^{(i,n+1-i)} \qquad  \forall 0<i<n $ resp. 

$ \qquad 0<i<n+1 $

is called $S$-compatible.

\end{Definition}

\begin{Theorem}
\label{zentral}
Let there be an $n \in \N_0$  with

(i) $B_1(S) = 0 = B_2(S) , A_n(S) = 0 = A_{n+1}(S)$,

(ii) The only $S$-compatible endomorphism $E : V^{ \x n+2} \rightarrow V^{ \x n+2}$ is the trivial     endomorphism $E = 0$.

\bigskip

\noindent Then $ A_n(S) = 0 \qquad \forall n \in \N_0 $ (and, as above, $B_1(S) = 0 $ ), so $ \rho_S $ is invariant under Kirby moves and forms therefore a 3-manifold invariant.
   
\end{Theorem}   

\bigskip

\noindent {\em Proof.\/} The proof is founded onto two lemmata. 

\begin{Lemma} 
If there is no nonzero $S$-compatible endomorphism in the i-th tensor potence then there is also no nonzero $S$-compatible endomorphism for $ m>i $.
\end{Lemma}

\noindent {\em Proof of the Lemma.\/} The idea of the proof is: Instead of regarding a given $ v^m \times v^m $-matrix $K$ examine its 
$v^n \times v^n$-submatrices. 
The conditions ($\a$),($\b$),($\gamma$) for $K$ will generate analogous conditions for the submatrices.

\bigskip

\begin{figure}[H]
\epsfig{file=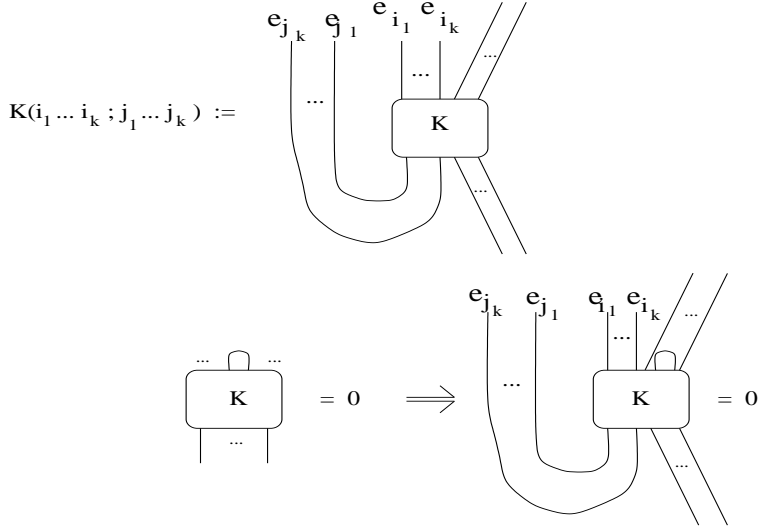,height=7cm,width=10cm,clip=,angle=0}
\caption{\footnotesize{The conditions go over to the submatrices}}
\label{untermatrizen}
\end{figure}

\noindent Define $K(i_1...i_{m-n-1};j_1...j_{m-n-1}) : V^{ \x n} \rightarrow V^{ \x n} $ as in fig. \ref{untermatrizen}: 

$ \langle K(i_1...i_{m-n-1};j_1...j_{m-n-1}) (v), w \rangle _n $

$:=\langle K(e_{i_1} \x ... \x e_{i_{m-n-1}} \x v , e_{j_1} \x ... \x e_{j_{m-n-1}} \x w \rangle _m $ . 

The positions of submatrices are as usual according to the v-adic system. 

As fig \ref{untermatrizen} illustrates, ($\a$),($\b$),($\gamma$) are transferred to the submatrices.

\noindent Therefore for all $ v^n \times v^n$-submatrices of $K$ the conditions in the definition hold, so they are equal to the zero matrix and therefore $K$, too \qed

\bigskip

\noindent It remains to show that the above differences $A_n$ are $S$-compatible in {\em each} tensor potence $n$. This will be guaranteed by the following Lemma.

\begin{Lemma}
The maps $A_n$ satisfy ($\a$),($\b$),($\gamma$) $ \forall n \in \N_0 $.
\end{Lemma}

\bigskip

\noindent {\em Proof of the Lemma.} (only to show for $m>n+2$)

\noindent The structure of the underlying induction is: 

$ ... \Rightarrow  A_{n-1}=0 \Rightarrow B_{n-1}=0 \Rightarrow A_n=0 \Rightarrow B_n=0 \Rightarrow ... $
  
\bigskip

\begin{figure}[ht]
\epsfig{file=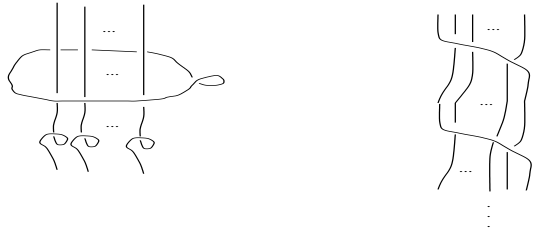,height=3cm,width=7cm,clip=,angle=0}
\caption{\footnotesize{Two tangles whose difference is $C_n$ }}
\label{Zopftwist}
\end{figure}

Instead of $A_n$, consider the map  
$ C_n := \rho(tw^{-1})^{ \x n} \o A_n $ (cf. fig. \ref{Zopftwist} ).  

\bigskip
The right tangle in the above figure we call {\em braid twist}. 

Now for the three conditions:

\begin{itemize}

\item{($\gamma$) is geometrically obvious.}

\item{ ($\a$) follows inductively from $A_{n-2} = 0$ (pull as above the two closed strands through the tangle; then the upper part is mapped to zero).} 

\item{For ($ \b $) cf. the following figures (! := to be deleted in the next step as closed component of a Fenn-Rourke move).}

\end{itemize}

\begin{figure}[H]
\epsfig{file=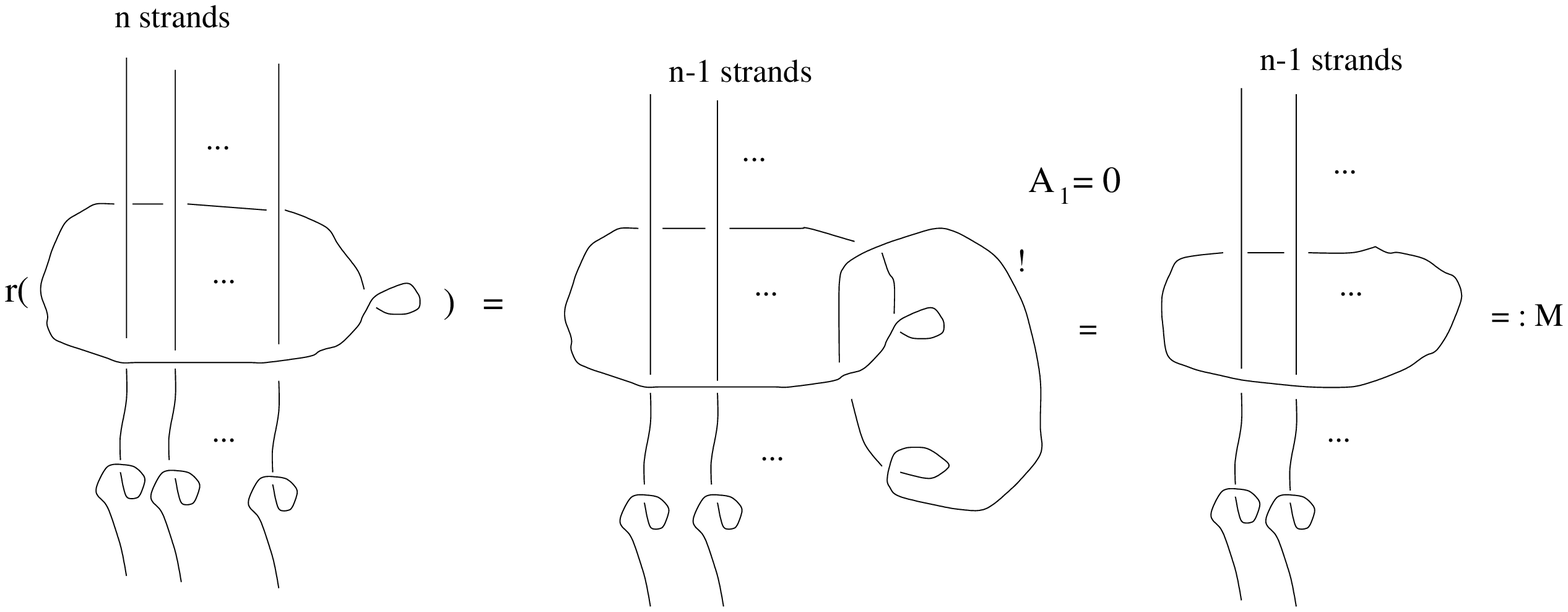,height=4cm,width=12cm,clip=,angle=0}
\caption{\footnotesize{The r-trace of the left tangle}}
\label{schlingling}
\end{figure}

\bigskip

\begin{figure}[H]
\epsfig{file=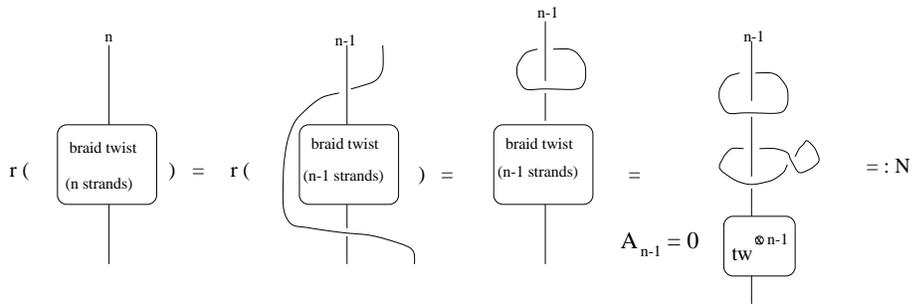,height=4cm,width=12cm,clip=,angle=0}
\caption{\footnotesize{The r-trace of the right tangle}}
\label{schlingrecht}
\end{figure}

\noindent In figure \ref{schlingling}, the last equation holds because of $A_1=0$, in figure \ref{schlingrecht},  the last equation holds because of $A_{n-1}=0$. 
Now we have to show $L := M-N = 0$. 
Consider $ L: V^{ \x n-1 } \rightarrow V^{ \x n-1 } $ satisfying ($\a$) because of the induction assumption for $n-2$, ($\gamma$) as above, for ($\beta$) apply again $r$ to both summands:

\begin{figure}[H]
\epsfig{file=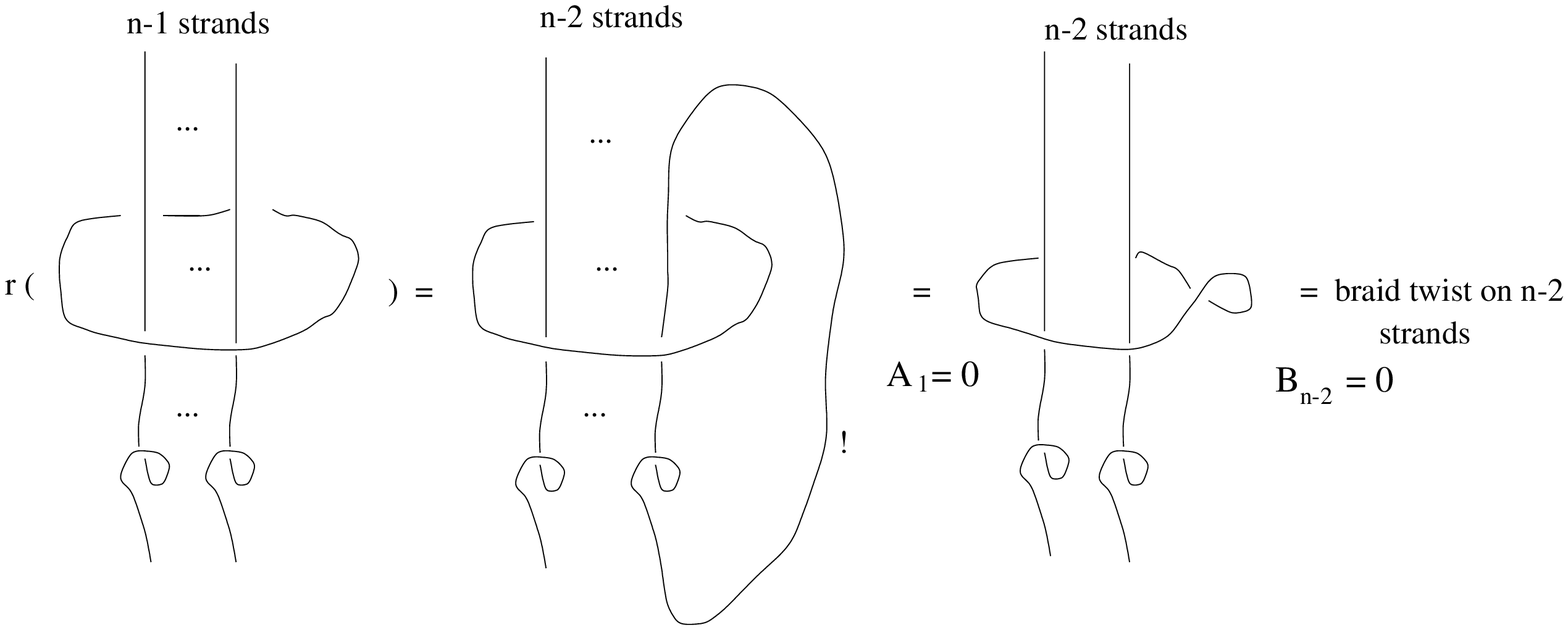,height=4cm,width=12cm,clip=,angle=0}
\caption{\footnotesize{The double r-trace of the left tangle}}
\label{schlingl2}
\end{figure}

\begin{figure}[H]
\epsfig{file=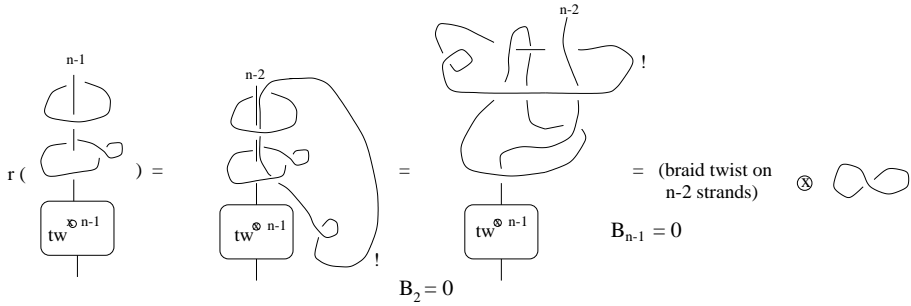,height=4cm,width=12cm,clip=,angle=0}
\caption{\footnotesize{The double r-trace of the right tangle}}
\label{schlingr2}
\end{figure}

\noindent In figure \ref{schlingl2}, the  last equation holds because of $A_1=0$, in figure \ref{schlingr2}, the last equation holds because of $B_{n-1}=0$, the one before because of $B_2=0$.
The proof for the $B_n$ is analogous with mirror images   \qed

\bigskip

\begin{Theorem}
\label{sym}

Conditions as in Theorem \ref{zentral}, but additionally $[S,T] =0 $ and instead of (ii)

(ii') The only map $ H : V^{ \x n+2 } \rightarrow V^{ \x n+2 } $ or 

$ H : V^{ \x n+3 } \rightarrow V^{ \x n+3 } $ with

$\qquad$ ($\a '$) $ b^{(i,n-i)} \o H = 0 $, $H$ symmetrical
 
$\qquad$ ($\b '$) $ l(H) = 0 $
 
$\qquad$ ($\gamma ' $)  H commutes with $ S^{(i,n-i)} $ and $ tw^{(i,n+1-i)} \qquad \forall 0<i<n $ 

$\qquad$ resp. $ 0<i<n+1 $

  are the trivial maps $H = 0$.

\noindent Then $ A_m (S) = 0 $ as above and $ \rho (S) $ generates a 3-manifold invariant.
\end{Theorem}

\noindent {\em Proof.} We have an analogon to the second lemma in the last proof: $ A_m, B_m $ are obviously T-symmetrical (because both summands are $s$-symmetrical). Hence we consider $ X := T_n \o A_n $. 
$X$ is symmetrical. The three conditions can be transferred to $X$:

\noindent ($\a '$) $ b^{(i,n-i)} \o H = 0 \Rightarrow  b^{(n-i-1,i+1)} \o (T_n \o H) = 0 $
 
\noindent ($\b '$) $ r(H) = 0 \Rightarrow  l(T_n \o H ) = 0 $ 

\noindent ($\gamma$) The commutating relations also can be transferred by index changes as in ($\a$). Here you need that $S$ and $T$ commute. 

\bigskip

\noindent There is also an analogon to the first lemma: Here we pass over from $d^n \times d^n$- and $d^{n+1} \times d^{n+1}$-matrices inductively to $d^{n+2} \times d^{n+2}$-matrices, i.e. if (ii) holds for $m, m+1$ then also for $m+2$.
 Again, we subdivide the matrix in blocks satisfying all conditions like in the proof above (excepted the symmetry), so they are antisymmetrical (if such a submatrix, say $A$, is not, consider $(A + A^+)$ which is a symmetrical matrix and still satisfies all other conditions, in contradiction to assumption (ii) ). On the other hand, we can subdivide the matrix into $d \times d$-blocks and ttreat them as entries from the ring of $d \times d$-matrices forming a big $d^{n+1} \times d^{n+1}$-matrix. Again all conditions can be related to the $d \times d$-matrices because every excepted each last condition multiplicate the blocks as wholes. Therefore $K : V^{ \x n+2 } \rightarrow V^{ \x n+2 }$ is symmetrical in the $d \times d$-blocks, the $d^{n+1} \times d^{n+1}$-blocks are antisymmetrical in the numbers, and (what is equivalent to (ii) for n) the  $d^{n+1} \times d^{n+1}$-blocks are antisymmetrical with respect to the $ d \times d$-blocks. So in all $K$ is antisymmetrical. Together with the condition that $K$ is symmetrical this means $K = 0$ \qed

\section{An Example}

\begin{Theorem}

Let $S$ be an irreducible S-Matrix with

$tr (S^{ \pm 1} ) = 1$, 

$A_1(S) = 0 = B_1(S)$. 

\noindent Then $ \rho_S $ is invariant under Kirby moves and induces therefore a 3-manifold invariant.

\end{Theorem}

\noindent {\em Proof.} The conditions on the trace are equivalent to $ A_0(S) = 0 = B_0 (S) $.

\noindent Apply Theorem \ref{zentral}:

\begin{itemize}

\item{ $E$ does not have full rang (because of $ (\alpha ) $) }

\item{ If $ v \notin im(E)$ then also $ Sv, S^{-1}v \notin im(E) $ (because of $ (\gamma )$ )}

\item{Because of irreducibility, in the rational normal form $S$ has got only one block; at least one corresponding basis vector is not contained in the image of $E$, so no one is

$\Rightarrow   im(E) = 0 $           \qed}

\end{itemize}

\section{Perspectives and Open Questions}

The presented theorems make an algorithm possible that checks solutions of the YBE on the criteria mentionned above. But the question of {\em existence} of such a solution is still open. For dim(V)=2 we have a definitely negative answer (not due two the special construction but for the general reasons mentionned above). For dim(V)=3 already the classification of the solutions of the YBE is a nontrivial problem of computer algebra (compare ~\cite{do}).

\noindent But {\em if} there is such an endomorphism the induced invariant $ \i $ will distinguish infinitely many 3-manifolds because

$ \displaystyle{ \i (  (\S^1 \times \S^2)^{ \# n} ) = \rho_S ( \bigsqcup_{i=1}^n 0 ) = (dim (V) )^n }$.


\begin{thebibliography}{99}



\begin{small}

\bibitem{do}
Christian Doll: Verschlingungsinvarianten aus neuen L\"osungen der Yang-Baxter-Gleichung. 
Diplomarbeit, Universit\"at Bonn, 1998

\bibitem{fr}
Roger Fenn/Colin Rourke: On Kirby's calculus of links, in: Topology 18 (1979), S. 1-15

\bibitem{fy}
Peter J. Freyd, David N. Yetter: Braided Compact Closed Categories with Applications to Low Dimensional Topology, in: Advances in Mathematics 77 (1989), S. 156-182

\bibitem{ki}
R.Kirby, A calculus for framed links in $ \S^3$, in: Invent. Math. 45 (1978), S. 35-56

\bibitem{nl}
Ning Lu: A simple proof of the fundamental theorem of Kirby calculus of links, in: Transactions of the American Mathematical Society 331 (1992), S. 143-156

\bibitem{m}
Olaf M\"uller: Geometrische Invarianten aus monoidalen Darstellungen. Diplomarbeit, Universit\"at Bonn, 2000

\bibitem{ps}
V.V. Prasolov/A.B. Sossinsky: Knots, Links, Braids and 3-Manifolds. Translations of Mathematical Monographs, American Mathematical Society 1991

\bibitem{t2}
V.G. Turaev: Operator invariants of tangles, and R-matrices, 
in: Math. USSR Izvestiya 35 (1990), No. 2, S. 411-444

\bibitem{t3}
V.G. Turaev: Quantum Invariants of Knots and 3-Manifolds. Walter de Gruyter, Berlin, New York 1994

\bibitem{we}
Hans Wenzl: Subfactors and Invariants of 3-Manifolds. Operator Algebras, Mathematical Physics, and Low Dimensional Topology (1993)

\end{small}

\end{thebibliography}
\end{document}